\documentclass[12pt,reqno]{amsart}
\usepackage{fullpage}
\usepackage{lineno,hyperref}
\usepackage{tikz}
\usetikzlibrary{matrix,arrows,decorations.pathmorphing}
\usepackage{tikz-cd}

\usepackage{amsmath}
\usepackage{epsfig}
\usepackage{amsmath,amsfonts,amsthm,amscd, amssymb}
\usepackage{latexsym}
\usepackage{enumerate}
\usepackage{mypack}
\usepackage[parfill]{parskip}

\title{Colimits of abelian groups}
\date{}
\address{Department of Mathematics, University of Western Ontario, London, Ontario, Canada N6A 5B7}\email{cokay@uwo.ca}
\author{C\.{I}han Okay}
\begin{document}
\maketitle

\begin{abstract}
In this paper we study the colimit $N_2(G)$ of abelian subgroups of a discrete group $G$. This group is the fundamental group of a subspace $B(2,G)$ of the classifying space $BG$ as introduced in  \cite{adem1}.  We describe  $N_2(G)$ for certain groups, and apply our results to study the homotopy type of the space $B(2,G)$. We give a list of classes of groups for which $B(2,G)$ is \underline{not} an Eilenberg--Maclane space of type $K(\pi,1)$.
\end{abstract}

\section{Introduction}
It is well known that the classifying space $BG$ of a discrete group $G$ can be described as a simplicial set where the $n$--simplices are ordered $n$--tuples of  elements in the group. It is possible to introduce certain simplicial subsets by putting conditions on the ordered $n$--tuples which respect the simplicial structure. An example of such a  construction as given in \cite{adem1} is to consider the simplicial set $B(q,G)$ whose $n$--simplices consist of $n$--tuples $(g_1,g_2,\cdots,g_n)$ of elements in $G$ such that these elements generate a subgroup $\Span{g_1,g_2,\cdots,g_n}\subset G$ of nilpotency class less than $q$. For example $B(2,G)$ consists of $n$--tuples of pairwise commuting elements. Equivalently, we can say that the set of $n$--simplices is the space of homomorphisms $\Hom(\ZZ^n,G)\subset G^n$. This definition also works more generally for topological groups where the space of homomorphisms inherits a  topology from the direct product $G^n$. For a recent survey on spaces of homomorphisms and related constructions such as $B(q,G)$ see \cite{cohen}. A key property of $BG$ is that it classifies principal $G$--bundles; there is a bijection between the set of isomorphism classes of principal $G$--bundles over a CW--complex $X$ and the set of homotopy classes of maps $X\rightarrow BG$. It is shown in \cite{adem2} that $B(q, G)$ is a classifying space for principal $G$--bundles of certain type, so called principal $G$--bundles of transitional nilpotency class less than $q$.
It is a natural problem to study the homotopy type of $B(q,G)$. It is \underline{not} always true that for a discrete group $G$ the space $B(q,G)$ is an Eilenberg--Maclane space of type $K(\pi,1)$, as opposed to the case of the usual classifying space. It was conjectured in \cite[page 15]{adem1} that for finite groups $B(q,G)$ is a $K(\pi,1)$ space, a natural expectation as a consequence of the types of groups studied there. Extraspecial $2$--groups appeared as the first counter examples in \cite{okay}. Our approach is purely algebraic, except in the last section we apply our results to $B(q,G)$ to study its homotopical properties.

The connection to algebra is via the fundamental group. In \cite{adem1} it is shown that the fundamental group of $B(q,G)$ is isomorphic to the colimit $N_q(G)$ of the nilpotent subgroups $N\subset G$  of class less than $q$. It turns out that some of the basic homotopy theoretic properties of $B(q,G)$ are determined by the group $N_q(G)$. Throughout the paper we focus on the case $q=2$. Then $N_2(G)$ has the following presentation
$$
\Span{(g),\;g\in G|\;(g)(h)=(gh)\text{ if } [g,h]=1}
$$  
where $[g,h]=ghg^{-1}h^{-1}$. The natural map $\epsilon:N_2(G)\rightarrow G$ induced by the inclusion $B(2,G)\subset BG$ can be described by the assignment $(g)\mapsto g$. Let $D_2(G)$ denote the kernel of the homomorphism $G\times G\rightarrow G/[G,G]$ defined by the multiplication map $(x,y)\mapsto xy[G,G]$. There is also a natural map $\bar{\epsilon}:N_2(G)\rightarrow D_2(G)$ defined by $(g)\mapsto (g,g^{-1})$, which factors $\epsilon$. The main result of the paper essentially gives a group theoretic condition which implies that the map $\bar{\epsilon}$ is an isomorphism. This condition is given by a set of non--identity elements $\set{g_i}_{i=1}^{2r}$, which we call a symplectic sequence, satisfying the commutation rules:
\begin{align*}
&[g_i,g_{i+r}]=[g_j,g_{j+r}]  \text{ for all  } 1\leq i,j\leq r, \\
&[g_i,g_j]=1 \text{ for any other pair.}
\end{align*} 
The sequence is called non--trivial if $[g_i,g_{i+r}]\neq 1$.
\begin{thm}\label{thm1}Let $\set{g_i}_{i=1}^{2r}$ be a non--trivial symplectic sequence in $G$ for some $r\geq 2$, and $S$ denote the subgroup generated by $\set{g_i}$. Then the natural map $N_2(S)\rightarrow D_2(S)$ is an isomorphism. Moreover $N_2(S)\rightarrow N_2(G)$ is injective, and its image is the subgroup generated by $\set{(g_i)}$.
\end{thm}
When $G$ is finite this theorem gives an explicit group theoretic condition which implies that the space $B(2,G)$ is \underline{not} an Eilenberg--Maclane space of type $K(\pi,1)$. 
\begin{thm}\label{thm2}
Suppose that $G$ is a finite group which has a non--trivial symplectic sequence $\set{g_i}_{i=1}^{2r}$ for some $r\geq 2$. Then $B(2,G)$ is \underline{not} a $K(\pi,1)$ space. 
\end{thm}
Using this theorem we extend the list of groups $G$ for which $B(2,G)$ is not a $K(\pi,1)$ space to the following classes of groups: Extraspecial $p$--groups of rank $\geq 4$, general linear groups $\GL_n(\FF_q)$, $n\geq 4$, over a field of characteristic $p$, and symmetric groups $\Sigma_{k}$ on $k$ letters where $k\geq 2^4$.

The organization of the paper is as follows. In Section \ref{section:pre} we describe the basic properties of the groups $N_q(G)$. In the rest of the paper we mostly restrict to the case $q=2$. Symplectic sequences are defined, and  Theorem \ref{thm1} is proved in Section \ref{section:symplectic}. Applications of the algebraic results, in particular Theorem \ref{thm2}, and examples are given in Section \ref{section:examples}.
\\

\textbf{Acknowledgments.} The author would like to thank the referee for their helpful comments.

\section{\label{section:pre}Preliminaries}

The descending central series of a group $G$ is the normal series
$$
1\subset\cdots\subset \Gamma ^{q+1} (G) \subset  \Gamma ^{q} (G) \subset \cdots \subset  \Gamma ^{2} (G) \subset \Gamma ^{1} (G)=G
$$
defined inductively $\Gamma ^1(G)=G$, $ \Gamma ^{q+1} (G)= [ \Gamma ^{q} (G) ,G]$.  A group is  nilpotent of class less than $q$ if  $\Gamma^{q}(G)=1$. The collection $\nN(q,G)$ of subgroups of class less than $q$ is a partially ordered set under inclusion. The colimit of the groups in this poset in the category of groups is denoted by 
$$N_q(G)= \colim{\nN(q,G)}N .$$ The colimit comes with a collection of homomorphisms $\eta_N:N\rightarrow N_q(G)$ such that the diagram 
$$
\begin{tikzcd}
N_q(G) \arrow{r}{\epsilon} & G \\
N \arrow[hook]{u}{\eta_N}\arrow[hook]{ru}
\end{tikzcd}
$$
commutes. In particular this implies that $\eta_N$ is an inclusion for every $N$.
There is a natural map
$
\epsilon:N_q(G)\rightarrow G
$
induced by the inclusions $N\rightarrow G$. This map is surjective since all the cyclic subgroups $\Span{g}$ are contained in $\nN(q,G)$, where $g\in G$.  We denote the image of $g$ under $\eta_{\Span{g}}$ simply by $\eta(g)$ or $(g)$. The construction $G\mapsto N_q(G)$ is natural and defines an endofunctor $N_q$ of the category of groups $\catGrp$. 

The group $N_q(G)$ can be constructed as a quotient of the free group $F(G)$ generated on the set  $\set{(g)|\;g\in G}$. Let  $R_q$ denote the normal closure of the subgroup of $F(G)$ generated by the products $(gh)^{-1}(g)(h)$   such that the subgroup $\Span{g,h}\subset G$  is of class less than $q$. The natural map  $(g)\mapsto \eta(g)$ establishes an isomorphism between $F(G)/R_q$ and $N_q(G)$. Its inverse is induced by the maps $N\rightarrow F(G)/R_q$. Therefore we obtain a presentation of the colimit 
\begin{equation}\label{presentation}
 N_q(G)=\Span{(g),\;g\in G|\;\;(g)(h)=(gh)\;\;\text{ if }\;\;\Gamma^q(\Span{g,h})=1 }.
\end{equation} 
When $q=2$ the relations remember the multiplication between the commuting elements, as $q$ gets larger more relations are added which correspond to higher nilpotence information. In particular we have a sequence of surjections
 $$
N_2(G)\twoheadrightarrow N_3(G)\twoheadrightarrow \cdots \twoheadrightarrow N_{q}(G)\twoheadrightarrow N_{q+1}(G)\twoheadrightarrow \cdots \twoheadrightarrow G .
 $$ 
Alternatively these groups can be described by their universal property: There is a natural bijection
$$
\catGrp(\colim{\nN(q,G)}N,H)\cong \ilim{\nN(q,G)}{}\catGrp(N,H).
$$
This means that a homomorphism $N_q(G)\rightarrow H$ is a set map $\phi:G\rightarrow H$ which restricts to a group homomorphism on the members of $\nN(q,G)$, we call such a map a \textit{$\nilq$--map}. Equivalently it is a set map satisfying $f(g_0)f(g_1)=f(g_0g_1)$ whenever $\Gamma^q(\Span{g_0,g_1})=1$. Therefore enlarging the morphisms of the category of groups by allowing such maps is equivalent to studying the image of the functor $N_q$. Define a category $\catGrp_q$ whose objects are groups and morphisms are $\nilq$--maps. Note that $N_q$ extends to this category and is the left adjoint of the inclusion functor
$$
N_q:\catGrp_q\rightleftarrows \catGrp:\iota_q
$$
with the unit $\eta:G\rightarrow N_q(G)$ defined by $\eta(g)=\eta_{\Span{g}}(g)$ and the counit $\epsilon:N_q(G)\rightarrow G$. Any  $\nilq$--map  $\phi:G\rightarrow H$ induces a group homomorphism $\epsilon N_q(\phi):N_q(G)\rightarrow H$. Note also that any $\nilr$--map is a $\nilq$--map for $r\geq q$. There is a sequence of inclusions of categories
$$
\catGrp\subset\cdots\subset \catGrp_{q+1} \subset  \catGrp_q \subset \cdots \subset  \catGrp_3 \subset \catGrp_2.
$$
 
Note that if a group $G$ is nilpotent of class less than $q$ then every $\nil_q$--map $G\rightarrow H$ is by definition a group homomorphism. Conversely, it is tempting to characterize groups $G$ for which every $\nil_q$--map $G\rightarrow H$ is a group homomorphism. We state the following conjecture as a new characterization of nilpotent groups:
\Conj{\label{q}A group $G$ is nilpotent of class less than $q$ if and only if the natural map $\epsilon:N_q(G)\rightarrow G $ is an isomorphism.}
We prove the case  $q=2$ which is special since there is a distinguished class of $\nil_2$--maps $\omega_n:N_2(G)\rightarrow N_2(G)$ for each integer $n\in \ZZ$ defined  by raising a generator  to the $n$--th power $(g)\mapsto (g^n)$. The map $\omega_{-1}$  is an automorphism  of $N_2(G)$ induced by the inversion map $i:G\rightarrow G$ which sends an element to its inverse.  
\Pro{\label{q=2}A group $G$ is abelian if and only if the natural map  $\epsilon:N_2(G)\rightarrow G$ is an isomorphism.}
\begin{proof}The inversion map $i$ induces a group homomorphism $\omega_{-1}$ on $N_2(G)$ hence by the commutative diagram of $\nil_2$--maps
$$
\begin{tikzcd}
G \arrow{r}{i} \arrow{d}{\eta}
&G \\
N_2(G) \arrow{r}{\omega_{-1}} &N_2(G)\arrow{u}{\epsilon}
\end{tikzcd}
$$
it induces a group homomorphism on $G$.
\end{proof}
 
To study the $q=2$ case in more detail we introduce the group $D_2(G)$. The group $D_2(G)$ is defined to be the kernel of the multiplication map
$$
G\times G\rightarrow G/[G,G],\;\;\;\;(x,y)\mapsto xy[G,G].
$$
Observe that the natural map $\epsilon:N_2(G)\rightarrow G$ which sends a generator $\eta(g)$ to $g$ factors as
$$
\begin{tikzcd}
N_2(G) \arrow{rd}{\epsilon} \arrow{d}{\bar{\epsilon}} \\
D_2(G) \arrow{r}{\pi_1} & G 
\end{tikzcd}
$$
where $\bar{\epsilon}$ sends $\eta(g)$  to $(g,g^{-1})$, and $\pi_1$ is the projection onto the first coordinate. 

\Pro{\label{surjective} The map $\bar{\epsilon}$ is surjective, and $\epsilon$ is the composition $\pi_1\bar{\epsilon}$.}
\Proof{The kernel of the projection map $\pi_1:D_2(G)\rightarrow G$ is the subgroup $1\times [G,G]$.
The surjectivity of $\bar{\epsilon}$ follows from the equation
$$
((gh)^{-1},gh)(g,g^{-1})(h,h^{-1})=(1,[g,h]),
$$
which implies that $D_2(G)$ is generated by the pairs of the form $(g,g^{-1})$. Second part of the statement is clear from the definition of the maps.
}

\section{\label{section:symplectic}Symplectic sequences}

In this section we study the group $N_2(G)$ by using  the presentation \ref{presentation} given in \S \ref{section:pre}:
$$
\Span{(g)|\;(g)(h)=(gh)\text{ if }[g,h]=1}.
$$
The main result of this section is a group theoretic condition which implies that the map  $\bar{\epsilon}:N_2(G)\rightarrow D_2(G)$ is  an isomorphism. Next definition is essential in our study.
\Def{\rm{ A sequence of  elements $\set{g_i}_{i=1}^{2r}$ in $G$ is called a \textit{symplectic sequence} if the following conditions are satisfied
\begin{enumerate}
\item $c=[g_i,g_{i+r}]=[g_j,g_{j+r}]$ for all $1\leq i,j\leq r$, 
\item  $[g_i,g_j]=1$ for all $1\leq i,j\leq 2r$ and $|i-j|\neq r$.
\end{enumerate}
The sequence is called \textit{non--trivial} if  $c\neq 1$, and \textit{trivial} if $c=1$.
 }}
Note that the element $c=[g_i,g_{i+r}]$ commutes with all the other $g_j$ for all $j$. 
The important fact about  symplectic sequences is that they are preserved under $\nil_2$--maps. This will be a special case of the following computation.
\Lem{\label{computation}
Let $\set{g_i}_{i=1}^{2r}$ be a symplectic sequence in $G$, and $\phi:G\rightarrow H$ a $\nil_2$--map. If $r\geq 2$ then for any positive integer $a,b,c,d$ the equation
$$
[\phi(g_i^ag_{i+r}^b),\phi(g_i^cg_{i+r}^d)]=[\phi(g_j),\phi(g_j)]^{ad-bc},\;\; 1\leq i,j\leq r,
$$
holds in $H$.
}
\Proof{For simplicity of notation set  $t_i=g_i^ag_{i+r}^b$, $t_{i+r}=g_i^cg_{i+r}^d$, and $n=ad-bc$.
We will use the following observation
$$
\begin{array}{lll} 
[g_j^{-n}t_{i+r},t_i^{-1}g_{j+r}]&=&[g_j^{-n},g_{j+r}][t_{i+r},t_{i}^{-1}]\\
&=&[g_j^n,g_{j+r}]^{-1}[t_{i},t_{i+r}]=1.
\end{array}
$$
Note that the commutators above lie in the subgroup $\Span{c}\subset G$ where $c$ commutes with all $g_i$. For all $ j\neq i$ the following  holds  
$$
\begin{array}{lll}
[\phi(t_i),\phi(t_{i+r})]&=& \phi(t_i)\phi(g_j^ng_j^{-n})\phi(t_{i+r})\phi(t_i)^{-1}\phi(g_{j+r}g_{j+r}^{-1})\phi(t_{i+r})^{-1}\\
&=& \phi(g_j^n)\phi(t_i)\phi(g_j^{-n}t_{i+r})\phi(t_i^{-1}g_{j+r})\phi(t_{i+r})^{-1}\phi(g_{j+r}^{-1})\\
&=&\phi(g_j^n)\phi(t_i)\phi(t_i^{-1}g_{j+r})\phi(g_j^{-n}t_{i+r})\phi(t_{i+r})^{-1}\phi(g_{j+r}^{-1})\\
&=&\phi(g_j^n)\phi(g_{j+r})\phi(g_j^n)^{-1}\phi(g_{j+r})^{-1}\\
&=&[\phi(g_j^n),\phi(g_{j+r})]\\
&=&[\phi(g_j),\phi(g_{j+r})]^n
\end{array}
$$
where we also make use of the identity $\phi(g)^{-1}=\phi(g^{-1})$. 
}
Taking $(a,b)=(1,0)$ and $(c,d)=(0,1)$ we obtain the following.
\Cor{\label{preserve} 
Suppose that $\set{g_i}_{i=1}^{2r}$, $r\geq 2$, is a symplectic sequence in $G$, and $\phi:G\rightarrow H$ is a $\nil_2$--map. Then 
$\set{\phi(g_i)}_{i=1}^{2r}$ is a symplectic sequence in $H$.
}
For the rest of this section  assume that $G$ has a non--trivial symplectic sequence $\set{g_i}_{i=1}^{2r}$ where $r\geq 2$. Let $S\subset G$ denote the subgroup generated by the elements $\set{g_i}_{i=1}^{2r}$. We make some observations about the structure of $S$. The commutator of $S$ is generated by $c=[g_i,g_{i+r}]$, and is contained in the center $Z(S)$. This implies that the commutator is bilinear on $S$: $[xy,z]=[x,z][y,z]$. Recall that we write $(x)$ to denote the image of $x$ under the canonical $\nil_2$--map $\eta:S\rightarrow N_2(S)$. According to the presentation \ref{presentation} we are allowed to write
$$
(x)(y)=(xy) \text{ if } [x,y]=1.
$$
In particular $(x^m)=(x)^m$ for any $m$. Moreover the commutator subgroup of $N_2(S)$ can be computed by applying Lemma \ref{computation} to the $\nil_2$--map $\eta$:
\begin{equation}\label{bilinear}
[(g_i^ag_{i+r}^b),(g_i^cg_{i+r}^d)]=[(g_j),(g_{j+r})]^{ad-bc}.
\end{equation}
Therefore  the commutator subgroup  is a cyclic group generated by $[(g_i),(g_{i+r})]$, and it is contained in the center of the group. To understand the structure of $N_2(S)$ we prove some relations between the elements $$k_i=(g_ig_{i+r})^{-1}(g_i)(g_{i+r}).$$ By equation \ref{bilinear} these elements are central in $N_2(S)$. 

\Lem{\label{computation2}For $1\leq i\neq j \leq r$  the relation $k_i=k_j$ holds in $N_2(S)$.
}
\Proof{Using the relations in $N_2(S)$ we have
$$
\begin{array}{lll}
k_ik_j^{-1}&=&(g_ig_{i+r})^{-1}(g_i)(g_{i+r})(g_{j+r})^{-1}(g_j)^{-1}(g_jg_{j+r}) \\
&=& (g_ig_{i+r})^{-1}(g_ig_{j+r}^{-1})(g_{i+r}g_j^{-1})(g_jg_{j+r}) \\
&=& (g_ig_{i+r})^{-1}(g_ig_{j+r}^{-1}g_{i+r}g_j^{-1})(g_jg_{j+r})\\
&=&(g_{i+r}^{-1}g_{i}^{-1}g_ig_{j+r}^{-1}g_{i+r}g_j^{-1})(g_jg_{j+r})\\
&=&(g_{j+r}^{-1}g_j^{-1})(g_jg_{j+r})=1.
\end{array}
$$
} 
Therefore we can simply omit the subscripts, and write $k=k_i$. 

\Lem{\label{power}
For $m\geq 1$ the  relation
$$
k^m=(g_i^{m}g_{i+r})^{-1}(g_i^m)(g_{i+r})
$$
holds in $N_2(S)$.
}
\Proof{
This is proved by induction on $m$. For $m=1$ the relation holds by definition. The result will follow from the computation
$$
\begin{array}{lll}
(g_i^{a}g_{i+r})^{-1}(g_i)(g_{i}^{a-1}g_{i+r})&=&(g_{j+r})(g_i^{a}g_{i+r})^{-1}(g_{j+r}^{-1}g_i)(g_i^{a-1}g_{i+r}g_{j}^{-1})(g_j)\\
&=&(g_{j+r})(g_i^{a}g_{i+r})^{-1}(g_{j+r}^{-1}g_i^{a}g_{i+r}g_{j}^{-1})(g_j)\\
&=&(g_jg_{j+r})^{-1}(g_j)(g_{j+r})  =k,
\end{array}
$$
where we used $[(g_{j+r}),(g_jg_{j+r})^{-1}(g_j)]=1$.
Now assuming the statement holds for $m=a-1$ we can write
$$
\begin{array}{lll}
(g_i^{a}g_{i+r})&=&k^{-1}(g_i)(g_i^{a-1}g_{i+r})\\
&=&k^{-1}(g_i)k^{1-a}(g_i^{a-1})(g_{i+r})\\
&=& k^{-a}(g_i^a)(g_{i+r}).
\end{array}
$$
}

\Lem{\label{merge}
The following relation holds in $N_2(S)$:
$$
(x)(y)(xy)^{-1}=k^\alpha\; \text{ for all } x,y\in \Span{g_i,g_{i+r}},
$$
 where $\alpha$ is such that $[x,y]=c^\alpha$.
}
\Proof{
For $j\neq i$ we have
$$
\begin{array}{lll}
(x)(y)&=&(g_{j+r})(g_{j+r}^{-1}x)(yg_j^{-\alpha})(g_j)^{\alpha}\\
&=&(g_{j+r})(g_{j+r}^{-1}xyg_j^{-\alpha})(g_j)^{\alpha} \\
&=&(xy)(g_{j+r})(g_{j+r}^{-1}g_j^{-\alpha})(g_j)^{\alpha} \\
&=&(xy)(g_{j+r}^{-1}g_j^{-\alpha})(g_j)^{\alpha} (g_{j+r})\\
&=&(xy)(g_j^\alpha g_{j+r})^{-1}(g_j^{\alpha})(g_{j+r})\\
&=&(xy) k^\alpha.
\end{array}
$$
In the last step we used Lemma \ref{power}.
}
Observe that any element $x$ in $S$  can be written as a product $c^kx_1x_2\cdots x_r$ for some $k$ where $x_i$ is of the form $g_i^{a_i}g_{i+r}^{b_{i+r}}$ for some $a_i, b_{i+r}$. Moreover, $[x_i,x_j]=1$ for all $i,j$. Therefore  $(x)=(c^k)(x_1)(x_2)\cdots (x_{r})$ since $[(x_i),(x_j)]=1$ in $N_2(S)$.
\Lem{\label{kernel}
The group $N_2(S)$ sits in a central extension
$$
1\rightarrow \Span{k}\rightarrow N_2(S)\stackrel{\epsilon}{\rightarrow} S\rightarrow 1.
$$
}
\Proof{
The kernel of $\epsilon$ is the normal closure of the subgroup generated by $(xy)^{-1}(x)(y)$ where $x,y\in S$. We can write $x=c^k\prod_{i=1}^r x_i $ and $y=c^l\prod_{i=1}^r y_i $ for some $k,l$, and $x_i,y_i\in \Span{g_{i},g_{i+r}}$.  Then we compute
$$
\begin{array}{lll}
(xy)^{-1}(x)(y) &=& (\prod x_iy_i)^{-1} (\prod x_i)(\prod y_i)(c)^{-k-l}(c)^k(c)^l\\
&=& \prod( x_iy_i)^{-1} \prod( x_i)\prod( y_i)\\
&=& \prod( x_iy_i)^{-1}( x_i)( y_i)\\
&=& \prod k^{\alpha_i}=k^{ \alpha_1+\alpha_2+\cdots+\alpha_r}
\end{array}
$$
where $\alpha_i$ is such that $[x_i,y_i]=c^{\alpha_i}$. We used Lemma \ref{merge} in the last step.
}

\Thm{\label{symplectic}
Let $\set{g_i}_{i=1}^{2r}$ be a non--trivial symplectic sequence in $G$ for some $r\geq 2$, and $S$ denote the subgroup generated by $\set{g_i}$. Then the natural map $N_2(S)\rightarrow D_2(S)$ is an isomorphism. Moreover $N_2(S)\rightarrow N_2(G)$ is injective, and its image is the subgroup generated by $\set{(g_i)}$.
}
\begin{proof} 
We want to show that the natural map 
$$
\bar{\epsilon}: N_2(S)\rightarrow D_2(S),\;\; (x)\mapsto (x,x^{-1}),
$$
is an isomorphism. This map is surjective by Proposition \ref{surjective}. We will show that it is also injective. Since $\bar{\epsilon}$ factors $\epsilon$ it suffices to show that $\ker(\epsilon)=\Span{k}$ maps injectively under $\bar{\epsilon}$.
 We have $\bar{\epsilon}(k)=(1,c)$. Assume $\bar{\epsilon}(k^m)=1$, that is $c^m=1$, for some $m$. 
In particular, this implies that $g_i^m$ lies in the center of $S$ for all $i$. 
 By Lemma \ref{power} we have 
 $$
\begin{array}{lll}
 k^m&=&(g_i^mg_{i+r})^{-1}(g_i^m)(g_{i+r})\\
 &=&(g_{i+r})^{-1}(g_i^m)^{-1}(g_i^m)(g_{i+r})=1
\end{array} 
$$
 since $[(g_i^m),(g_{i+r})]=1$. This proves the injectivity of $\bar{\epsilon}$. 

The second statement follows from the  diagram
$$
\begin{tikzcd}
N_2(S) \arrow{r} \arrow{d}{\cong}
&N_2(G) \arrow{d}\\
D_2(S) \arrow{r} &D_2(G)
\end{tikzcd}
$$
since the map $D_2(S)\rightarrow D_2(G)$ is injective.
\end{proof}
 
\section{\label{section:examples}Examples and applications}
In this section we list some of the consequences of the group theoretic results obtained in the previous sections. Our main result implies that there exists groups $G$ for which $B(2,G)$ is not a $K(\pi,1)$ space. On the geometric side, the space $B(2,G)$ is a  classifying space for principal $G$--bundles of transitional nilpotency class less than $2$.  A principal $G$--bundle  over a CW--complex $X$ is said to have transitional nilpotancy class less than $2$ if there exists a trivialization of the bundle by an open cover such that on intersections the transition functions commute when they are simultaneously defined.  Two such bundles $p_0$ and $p_1$ are called $2$--transitionally isomorphic if  there exits a principal $G$--bundle $p$ over $X\times [0,1]$ with transitional nilpotancy class less than $2$ whose restrictions $p|_{X\times 0}$ and $p|_{X\times 1}$  give $p_0$ and $p_1$, respectively \cite[Definition 5.3]{adem2}. One implication of having a non--trivial higher homotopy group $\pi_n(B(2,G))$ for some $n\geq 2$ is that there exists principal $G$--bundles over the $n$--sphere $\sS^n$ which are not $2$--transitionally isomorphic to the trivial principal $G$--bundle. 
Note that in the case of ordinary principal $G$--bundles every principal $G$--bundle over $\sS^n$ is isomorphic to the trivial principal  $G$--bundle when $n\geq 2$ since $\pi_n(BG)=0$ for all $n\geq 2$.

\subsection*{Homotopy type of $B(2,G)$}
The colimit $N_q(G)$ is isomorphic to the fundamental group of the geometric realization of a simplicial set $B(q,G)\subset BG$ introduced in \cite{adem1}. This simplicial set can be described as a colimit of the classifying spaces of nilpotent subgroups of class less than $q$ \cite[Theorem 4.3]{adem1}, also see \cite{okay} for further properties. Similar to the usual classifying space $BG$ the association $G\mapsto B(q,G)$ is functorial on the category of groups, moreover it extends to a functor  
$$
B(q,-):\catGrp_q\rightarrow \catS
$$
from the extended version of the category of groups as introduced in \S \ref{section:pre} to the category of simplicial sets. 

The inclusion map $B(q,G)\rightarrow BG$ defines a principal $G$--bundle  $E(q,G)\rightarrow B(q,G)$.  Looking at the long exact sequence of homotopy groups associated to the fibration in low degrees we obtain a short exact sequence  of groups 
\begin{equation}\label{short-exact}
1\rightarrow \pi_1(E(q,G))\rightarrow N_q(G)\rightarrow G\rightarrow 1,
\end{equation}
where we identify the fundamental group of $\pi_1(B(q,G))$ with the colimit $N_q(G)$.
The characterization in Proposition \ref{q=2} of abelian groups can be translated into the following statement.

\Pro{A group $G$ is abelian  if and only if $E(2,G)$ is contractible.}
\Proof{If $G$ is abelian then by the description of $B(2,G)$ as a colimit of the classifying spaces $BA$ of abelian subgroups $A$, we have an isomorphism $B(2,G)\cong BG$. Therefore the space $E(2,G)$ is isomorphic to the contractible  space $EG$. Conversely, if the space $E(2,G)$ is  contractible then the kernel of $N_2(G)\rightarrow G$ is trivial by \ref{short-exact}. By Proposition \ref{q=2} $G$ is abelian.
} 
More generally, one can ask when $B(q,G)$ has the homotopy type of a $K(\pi,1)$ space, a question raised in \cite{adem1}. There is a natural map $B(q,G)\rightarrow BN_q(G)$ as defined in \cite[Theorem 4.4]{adem1}.

\Pro{\label{torfree}Suppose that $G$ is a finite group. If the natural map $B(q,G)\rightarrow BN_q(G)$ is a homotopy equivalence then the kernel  of   $\epsilon:N_q(G)\rightarrow G$ is torsion free.}
\Proof{ Suppose that the map $B(q,G)\rightarrow BN_q(G)$ is a homotopy equivalence. This means that $B(q,G)$ is a $K(N_q(G),1)$ space. Since $E(q,G)$ is a finite cover of $B(q,G)$ this implies that $E(q,G)$ is a $K(\pi,1)$ space. The  space $E(q,G)$ is homotopy equivalent to the nerve of the poset of right cosets $\set{gN\subset G|\; N\in \nN_q(G)}$ of nilpotent subgroups of class less than $q$, ordered by inclusion  \cite[\S 3.3]{okay}. In particular, when $G$ is finite the nerve of this poset is finite dimensional hence  $E(q,G)$ has finite cohomological dimension. Therefore $E(q,G)$ has the homotopy type of a finite dimensional complex. Note that if a finite dimensional complex is a $K(\pi,1)$ space then its fundamental group is torsion free. This follows from the fact that the  cohomological dimension of $\pi$ is less than or equal to its geometric dimension, and $\pi$ is torsion free if its cohomological dimension is finite \cite[Chapter \rom{8}]{brown2}. Therefore  $\pi_1(E(q,G))$, which is isomorphic to the kernel of the map $N_q(G)\rightarrow G$, is torsion free.
}
Next we describe our main result on the homotopy type of $B(2,G)$. This theorem gives an algebraic condition on $G$ which implies that $B(2,G)$ is \underline{not} a $K(\pi,1)$ space. First examples of such groups are given in \cite{okay} using different methods.
\Thm{\label{eqtorsion}
Suppose that $G$ is a finite group which has a non--trivial symplectic sequence $\set{g_i}_{i=1}^{2r}$ for some $r\geq 2$. Then $B(2,G)$ is \underline{not} a $K(\pi,1)$ space. 
}
\begin{proof} Let $S\subset G$ denote the subgroup generated by $\set{g_i}$. By Theorem \ref{symplectic} there is a commutative diagram
$$
\begin{tikzcd}
N_2(S) \arrow[hook]{r} \arrow{d}{\epsilon_S} & N_2(G) \arrow{d}{\epsilon_G} \\
S\arrow[hook]{r} & G 
\end{tikzcd}
$$
and $\epsilon_S$ factors as $N_2(S)\cong D_2(S)\stackrel{\pi_1}{\rightarrow}S$. The kernel of $\pi_1$  is generated by the element $(1,c)$ which is of finite order since $G$ is finite. Hence there exists a torsion element in the kernel of $\epsilon_G$, and Proposition \ref{torfree} implies that $B(2,G)$ cannot be a $K(\pi,1)$ space.
\end{proof}

\subsection*{Examples}

We describe how to obtain a non--trivial symplectic sequence $\set{g_i}_{i=1}^{2r}$, $r\geq 2$, for the following groups so that Theorem \ref{eqtorsion} applies.
\begin{enumerate}
\item Extraspecial $p$--groups are central in the ideas developed in this paper, such as the definition of a symplectic sequence. An extraspecial $p$--group is a central extension
$$
Z(P)=\ZZ/p\lhd E_r \rightarrow (\ZZ/p)^{2r}
$$ 
where $Z(P)$ is also the commutator subgroup.
The quotient group has the structure of a symplectic vector space when regarded as a vector space over $\FF_p$ with the commutator map inducing a non--degenerate alternating bilinear form.
Assume that $r\geq 2$. The lift $\set{\tilde{e}_i}_{i=1}^{2r}$ of a symplectic basis $\set{e_i}_{i=1}^{2r}$ in the quotient group $(\ZZ/p)^{2r}$ gives a non--trivial symplectic sequence in $E_r$. Moreover, $E_r$ is generated by $\set{\tilde{e}_i}$ hence by Theorem \ref{symplectic} we have $N_2(E_r)\cong D_2(E_r)$. Therefore the kernel of $N_2(E_r)\rightarrow E_r$ is cyclic of order $p$.

\item The general linear group $\GL_n(\FF_q)$, $n\geq 4$, over a finite field $\FF_q$ of characteristic $p$    has a symplectic sequence given by the elementary matrices $$\set{g_1=E_{12},g_2=E_{13},g_3=E_{2n},g_4=E_{3n}}$$
 where $E_{ij}$ has $1$'s on the diagonal and  in the $(i,j)$-slot. 
 It follows easily from the commutation relations of elementary matrices that this sequence satisfies the required commutation relations of a non--trivial symplectic sequence.  
 
\item We can embed $\GL_n(\FF_p)$ inside the symmetric group $\Sigma_{p^n}$ 
$$
\iota:\GL_n(\FF_p)\rightarrow \Sigma_{p^n}
$$
by regarding linear transformations of an $n$--dimensional vector space $V$ over $\FF_p$ as the permutations of a set of cardinality $p^n$. When $n= 4$ the image of the symplectic sequence $\set{g_i}$ obtained from elementary matrices in $\GL_4(\FF_p)$ gives a non--trivial symplectic sequence in $\Sigma_{p^4}$. Therefore for $k\geq 2^4$ the natural inclusion $\Sigma_{2^4}\rightarrow \Sigma_{k}$ gives a non--trivial symplectic sequence in $\Sigma_k$.  One can conclude similarly for the alternating group $A_k$ by using the isomorphism    $\GL_4(\FF_2)\cong A_{8}$, and the natural inclusion $A_{8}\rightarrow A_{k}$ where $k\geq 8$. 
\end{enumerate}

\end{document}